\documentclass[a4paper]{article}
\usepackage[utf8x]{inputenc}
\usepackage[margin=1.0in]{geometry}

\usepackage{hyperref}
\title{On the relation between continuous functions in two different metric spaces}
\author{Adrian Fellhauer}
\usepackage{hyperref}
\usepackage{amssymb}
\usepackage{mathtools}
\begin{document}

\maketitle

\begin{abstract}
Let the metric space $\mathbb R^n \setminus \sim$ be the metric space of $n$-sized unordered tuples of real numbers. In the following, it will be shown that if a function $\varphi: \mathbb R^m \to \mathbb R^n \setminus \sim$ is continuous, then there is a continuous function $f: \mathbb R^m \to \mathbb R^n$ such that a natural embedding of $f$ into $\mathbb R^n \setminus \sim$ is equal to $\varphi$.

This theorem is wrong in the complex case. A counterexample is given in \cite{curgusmascioni06}.
\end{abstract}

\section{The metric space $\mathbb R^n \setminus \sim$}

Let $\mathfrak S_n$ be the set of permutations of size $n$. For every $\sigma \in \mathfrak S_n$, we define
\begin{equation*}
p_\sigma: \mathbb R^n \to \mathbb R^n  \text{, } p_\sigma((x_1, \ldots, x_n)) = (x_{\sigma(1)}, \ldots, x_{\sigma(n)})
\end{equation*}
We furthermore define the following set:
\begin{equation*}
\mathfrak P_n := \{p_\sigma : \sigma \in \mathfrak S_n\}
\end{equation*}
For $x, y \in \mathbb R^n$, we define the following equivalence relation:
\begin{equation*}
z \sim y :\Leftrightarrow \exists p_\sigma \in \mathfrak P_n : z = p_\sigma(y)
\end{equation*}

$\mathbb R^n \setminus \sim$ consists of the equivalence classes regarding this equivalence relation and is equipped with the following metric:
\begin{equation*}
d(\overline{y}, \overline{z}) := \min_{p_\sigma \in \mathfrak P_n} \|y - p_\sigma(z)\|_1
\end{equation*}
, where for $x \in \mathbb R^n$,
\begin{equation*}
\|x\|_1 := \sum_{k=1}^n |x_k|
\end{equation*}
is the 1-norm of $x$. $d$ has the following properties:
\begin{itemize}
\item It is well-defined, i. e. independent of the component's order
\item It is zero if and only if the two input elements are equal
\item It is symmetric, i. e. $d(\overline y, \overline z) = d(\overline z, \overline y)$
\item It fulfils the triangle inequality, i. e. $d(\overline y, \overline z) + d(\overline z, \overline x) \le d(\overline x, \overline y)$
\end{itemize}
In conclusion, we may say that $d$ is a metric, and $\mathbb R^n \setminus \sim$ is a metric space. The proof can be found in \cite[p.~391]{harrismartin87}.

\section{Three kinds of sets}

\subsection{Definitions}
Let $M_1, \ldots, M_i \subseteq \{1, \ldots, n\}$ such that $l \neq m \Rightarrow M_l \cap M_m = \emptyset$ and\\$\forall 1 \le l \le i : |M_l| \ge 2$. Then we define
\begin{equation*}
\mathcal X_{M_1, \ldots, M_i} := \{(z_1, \ldots z_n) \in \mathbb R^n | \forall 1\le j \le i : m, l \in M_j \Rightarrow z_m = z_l\}
\end{equation*}

\begin{equation*}
\mathcal Y_{M_1, \ldots, M_i} := \{p_\sigma \in \mathfrak P_n | \forall z \in \mathcal X_{M_1, \ldots, M_i} : p_\sigma(z) = z\}
\end{equation*}

\begin{equation*}
\mathcal X_{M_1, \ldots, M_i}^\epsilon := \{x \in \mathbb R^n | \min_{y \in \mathcal X_{M_1, \ldots, M_i}} \|x - y\| < \epsilon\}
\end{equation*}

\subsection{Lemma}
\begin{equation*}
\mathcal X_{M_1, \ldots, M_i}^\epsilon \subseteq \{x \in \mathbb R^n | \forall p_\sigma \in \mathcal Y_{M_1, \ldots, M_i} : \| p_\sigma(x) - x \| < 2 \epsilon \}
\end{equation*}

\subsubsection*{Proof}
Let $y \in \mathcal X_{M_1, \ldots, M_i}^\epsilon$. Then, by definition of $\mathcal X_{M_1, \ldots, M_i}^\epsilon$: $\exists x \in \mathcal X_{M_1, \ldots, M_i} : \| x - y\| < \epsilon$. But $p_\sigma(x) = x$, and due to commutativity of addition $\|x - y\| = \|p_\sigma(x) - p_\sigma(y)\|$, and therefore by the triangle inequality:
\begin{equation*}
\|y - p_\sigma(y)\| \le \|y- x\| + \|p_\sigma(y) - p_\sigma(x)\| < \epsilon + \epsilon = 2 \epsilon
\end{equation*}
$\Box$

\section{The set of non-descendingly ordered vectors}
\subsection{Definition}
We define the set of non-descendingly ordered vectors as follows:
\begin{equation*}
\mathbb R^n_o := \{(x_1, \ldots, x_n) \in \mathbb R^n : x_1 \le x_2 \le \cdots \le x_n\}
\end{equation*}

\subsection{Lemma}
$\mathbb R^n_o$ is closed. (Remark: This is the point where the theorem fails in the complex case. In fact, it has been shown that every set of complex numbers containing one element of each element in $\mathbb C^n \setminus \sim$ exactly once is not closed. See \cite[p.~12ff.]{curgusmascioni06})

\subsubsection*{Proof}
We show that the complement of $\mathbb R^n_o$ is open. Let $x = (x_1, \ldots, x_n) \notin \mathbb R^n_o$. Then, by defintion, we obtain
\begin{equation*}
\exists i, j \in \mathbb N : i < j \wedge x_i > x_j
\end{equation*}
Let $x_i - x_j =: c > 0$. Then we obtain $\forall y = (y_1, \ldots, y_n) \in B_{c/4}(x)$:
\begin{equation*}
|y_i - x_i| \le \|x - y\| < c/4 \text{ and } |y_j - x_j| \le \|x - y\| < c/4
\end{equation*}
, which is why
\begin{equation*}
y_i - y_j > c/2 \Rightarrow y \notin \mathbb R^n_o
\end{equation*}
$\Box$

\subsection{Lemma}
Let $B_\epsilon(y_0) \subset \text{int } \mathbb R^n_o$. Then it follows that
\begin{equation*}
\forall y \in B_\epsilon(y_0), p_\sigma \in \mathfrak P_n \setminus \{\text{Id}\} : p_\sigma(y) \notin \mathbb R^n_o
\end{equation*}

\subsubsection*{Proof}
Let $y \in B_\epsilon(y_0)$. We first notice that all components of $y$ are distinct, since if we assume the contrary, i. e. 
\begin{equation*}
y = (y_1, \ldots, y_i, \ldots, y_j, \ldots, y_n)
\end{equation*}
, where $y_i = y_j$, the sequence
\begin{equation*}
y_k = (y_1, \ldots, y_i + \frac{1}{k}, \ldots, y_j, \ldots, y_n)
\end{equation*}
converges to $y$ as $k \to \infty$, but lies outside $\mathbb R^n_o$, which is a contradiction to $B_\epsilon(y_0) \subset \text{int } \mathbb R^n_o$.

Let $p_\sigma \neq \text{Id}$. Then at least one component of $y$ changes position, meaning that the order is changed and therefore $p_\sigma(y) \notin \mathbb R^n_o$. $\Box$

\subsection{Lemma}
Let $y_0 \in \partial \mathbb R^n_o$. Then at least two components of $y_0$ are equal.

\subsubsection*{Proof}
Assume the contrary. Then $y_0$ would be of the form
\begin{equation*}
y_0 = (y^0_1, \ldots, y^0_n) \text{ with } y^0_1 < y^0_2 < \cdots < y^0_{n-1} < y^0_n
\end{equation*}
Let us choose
\begin{equation*}
0 < c := \min_{i \in \{1, \ldots, n-1\}} y^0_{i+1} - y^0_i
\end{equation*}
Then we obtain $\forall y = (y_1, \ldots, y_n) \in B_{c/4}(y_0), i \in \{1, \ldots, n\}$:
\begin{equation*}
|y_i - y^0_i| \le \|y - y_0\| < c/4
\end{equation*}
Therefore, we have
\begin{equation*}
\min_{i \in \{1, \ldots, n-1\}} y_{i+1} - y_i > c/2
\end{equation*}
, implying that $y$ is strictly ascendingly ordered. Therefore, $B_{c/4}(y_0) \subset \mathbb R^n_o$, and $y_0 \in \text{int } \mathbb R^n_o$, which contradicts the assumption. $\Box$

\subsection{Lemma}
Let $y_0 \in \partial \mathbb R^n_o$. If we order the sets $\mathcal X_{M_1, \ldots, M_i}$ with $y_0 \in \mathcal X_{M_1, \ldots, M_i}$ according to the partial order
\begin{equation*}
\mathcal X_{M_1, \ldots, M_i} \le \mathcal X_{W_1, \ldots, W_j} :\Leftrightarrow \mathcal X_{M_1, \ldots, M_i} \subseteq \mathcal X_{W_1, \ldots, W_j}
\end{equation*}
and choose a minimal element regarding this order $\mathcal X_{M_1, \ldots, M_i}$, we have
\begin{equation*}
\forall p_\sigma \in \mathfrak P_n \setminus \mathcal Y_{M_1, \ldots, M_i} : p_\sigma(y_0) \notin \mathbb R^n_o
\end{equation*}

\subsubsection*{Proof}
It is obvious that a minimum exists, since we are considering a finite and due to lemma 3.4 nonempty set.

Choose now $p_\sigma \in \mathfrak P_n \setminus \mathcal Y_{M_1, \ldots, M_i}$ arbitrarily. Then $y_0$ can not stay the same, since else by decomposition of $\sigma$ in disjoint cycles and iterated application of $p_\sigma$, we find that even more entries of $y_0$ must be equal (else $p_\sigma \in \mathcal Y_{M_1, \ldots, M_i}$). Therefore we would obtain a smaller set $\mathcal X_{W_1, \ldots, W_j}$ with $y_0$ in it, which is a contradiction to the assumption. But if $y_0$ changes, then $y_0$ is sorted differently after the application of $p_\sigma$, implying that $p_\sigma(y_0) \notin \mathbb R^n_o$. $\Box$

\section{Construction of a continuous function}

\subsection{Definition}
We define the following function:
\begin{equation*}
\psi: \mathbb R^n \setminus \sim \to \mathbb R^n_o \text{, } \psi(\overline x) = (x_1, \ldots, x_n) \text{ such that }(x_1, \ldots, x_n) \in \mathbb R^n_o \cap \overline x
\end{equation*}
$\psi$ is well-defined (i. e. $\forall \overline x \in \mathbb R^n \setminus \sim$ there is exactly one element in $\mathbb R^n_o \cap \overline x$) because otherwise there would be either no or two possibilities to sort a vector non-descendingly.

\subsection{Construction}
Let $\varphi: \mathbb R^m \to \mathbb R^n \setminus \sim$ be continuous. Then we construct $f: \mathbb R^m \to \mathbb R^n$ as follows:
\begin{equation*}
f: \mathbb R^m \to \mathbb R^n \text{, } f(x) := \psi(\varphi(x))
\end{equation*}
This function satisfies
\begin{equation*}
\forall x \in \mathbb R^m : \overline{f(x)} = \varphi(x)
\end{equation*}
This is why we say that a natural embedding of $f$ into $\mathbb R^n \setminus \sim$ is equal to $\varphi$.

\subsection{Proof of continuity of the constructed function}
For the sake of simplicity, in the following we will consider $\mathbb R^n$ equipped with the 1-norm. This can be done without losing generality because in finite dimensions, all norms are equivalent.\\\\
Let $x_0 \in \mathbb R^m$ be arbitrary, and let $y_0 := f(x_0)$. Let $\epsilon > 0$ be arbitrary. Since $\varphi$ is continuous, we may choose $\delta > 0$ such that
\begin{equation*}
\forall x \in B_\delta(x_0) : \varphi(x) \in B_\epsilon(\varphi(x_0))
\end{equation*}
Let $x \in B_\delta(x_0)$ be arbitrary. We consider two cases:\\\\
Case 1: $y_0 \in \text{int } \mathbb R^n_o$

In this case, we may choose $\epsilon$ small enough so that $B_\epsilon(y_0) \subset \mathbb R^n_o$ and then adjust $\delta$ accordingly. Lemma 3.3 implies that
\begin{equation}
\forall y \in B_\epsilon(y_0), p_\sigma \in \mathfrak P_n \setminus \{\text{Id}\} : p_\sigma(y) \notin \mathbb R^n_o
\end{equation}
Since $\varphi$ is continuous, we may choose $p_\sigma \in \mathfrak P_n$ such that $p_\sigma(f(x)) \in B_\epsilon(y_0)$. But since $f(x) = p_{\sigma^{-1}}(p_\sigma(f(x))) \in \mathbb R^n_o$ (1) implies $\sigma = \text{Id}$. Therefore we obtain $f(x) \in B_\epsilon(y_0)$, and since $\epsilon > 0$ and $x \in B_\epsilon(x_0)$ were arbitrary, continuity at $x_0$ is proven.\\\\
Case 2: $y_0 \notin \text{int } \mathbb R^n_o$

$y_0 \notin \text{int } \mathbb R^n_o$ means $y_0 \in \partial \mathbb R^n_o$, which is why $y_0 \in \mathcal X_{M_1, \ldots, M_i}$ for some $M_1, \ldots, M_i$ with $|M_l| \ge 2$, $ l \in \{1, \ldots, i\}$, $i \in \mathbb N$ (lemma 3.4). We order the sets $\mathcal X_{M_1, \ldots, M_i}$ with $y_0 \in \mathcal X_{M_1, \ldots, M_i}$ according to the order
\begin{equation*}
\mathcal X_{M_1, \ldots, M_i} \le \mathcal X_{W_1, \ldots, W_j} :\Leftrightarrow \mathcal X_{M_1, \ldots, M_i} \subseteq \mathcal X_{W_1, \ldots, W_j}
\end{equation*}
and choose one minimal element regarding this order $\mathcal X_{M_1, \ldots, M_i}$.

Since $\varphi$ is continuous, we may choose $p_\sigma \in \mathfrak P_n$ such that $p_\sigma(f(x)) \in B_\epsilon(y_0)$. Again $f(x) = p_{\sigma^{-1}}(p_\sigma(f(x)))$. If $p_{\sigma^{-1}} \in \mathcal Y_{M_1, \ldots, M_i}$, we obtain due to lemma 2.2 and the triangle inequality:
\begin{equation*}
\|f(x) - y_0\| \le \|p_{\sigma^{-1}}(p_\sigma(f(x))) - p_\sigma(f(x))\| + \|p_\sigma(f(x)) - y_0\| \le 3\epsilon
\end{equation*}
We lead $p_{\sigma^{-1}} \notin \mathcal Y_{M_1, \ldots, M_i}$ to a contradiction. Assume it were true. Then $p_{\sigma^{-1}}(y_0) \notin \mathbb R^n_o$ (lemma 3.5). Since $\mathbb R^n_o$ is closed (lemma 3.2) and since there are only finitely many $p_\sigma \in \mathfrak P_n \setminus \mathcal Y_{M_1, \ldots, M_i}$, we may choose $\eta_{y_0} > 0$ such that
\begin{equation*}
\forall p_\sigma \in \mathfrak P_n \setminus \mathcal Y_{M_1, \ldots, M_i} : B_{\eta_{y_0}}(p_\sigma(y_0)) \not\subset \mathbb R^n_o
\end{equation*}
If we now choose $\epsilon \le \eta_{y_0}$ and adjust $\delta$ accordingly, we find that for $x \in B_\delta(x_0)$, if $p_\sigma \in \mathfrak P_n$ such that $p_\sigma(f(x)) \in B_\epsilon(y_0)$ and $p_{\sigma^{-1}} \notin \mathcal Y_{M_1, \ldots, M_i}$, then
\begin{equation*}
\eta_{y_0} \ge \epsilon > \|f(x_0) - p_\sigma(f(x))\| = \|p_{\sigma^{-1}}(f(x_0)) - f(x)\|
\end{equation*}
and therefore $f(x) \notin \mathbb R^n_o$, which contradicts the definition of $f$ for small enough $\epsilon > 0$. Since $x_0$ was arbitrary and $\epsilon > 0$ arbitrarily small, we have proven that $f$ is continuous. $\Box$


\begin{thebibliography}{9}

\bibitem{curgusmascioni06}
  Branko Curgus and Vania Mascioni
,
  \emph{Roots and polynomials as homeomorphic spaces}.
  Expositiones Mathematicae, 
  Volume 42,
  Issue 1,
  February 2006,
  pp. 81-95.

\bibitem{harrismartin87}
  Gary Harris and Clyde Martin,
  \emph{The roots of a polynomial vary continuously as a function of the coefficients}.
  Proceedings of the American Mathematical Society
  Volume 100,
  Number 2,
  June 1987,
  pp. 390 - 392.

\end{thebibliography}
\end{document}